\documentclass[a4paper,11pt]{article}
\usepackage[utf8]{inputenc}
\usepackage{amsmath}
\usepackage{amsfonts}
\usepackage{amsthm}
\usepackage{multirow}
\author{Agnieszka Ka\l{}amajska $^{\rm a, b}$\thanks{Corresponding author. Email: A.Kalamajska@mimuw.edu.pl
\vspace{6pt}}
\ and Anna Maria Kosiorek$^{\rm a}$\\
\hspace{6pt}%$^{,}$
\vspace{-6pt}
\small{$^{a}${\em{Faculty of Mathematics, Informatics, and Mechanics, University of Warsaw, Poland}}};\\
\\
\small{$^{\rm b}${\em{Institute of Mathematics of the Polish Academy of Sciences at Warsaw, Poland}}}}

\date{}
\title{On maximum principles for radial solutions to nonlinear elliptic PDE's via Opial-type inequalities}

\newtheorem{tw}{Theorem}

\newtheorem{lem}{Lemma}
\theoremstyle{remark}
\newtheorem{re}{Remark}

\newtheorem{ex}{Example}

\renewcommand{\t}{\tau}

\newcommand{\la}{\lambda}

\newcommand{\R}{\mathbb{R}}

\newcommand{\h}{h(\cdot,\cdot,\cdot)}

\usepackage{enumitem}

\begin{document}
\date{\today}
\maketitle
\begin{abstract}
We consider degenerated nonlinear PDE of elliptic type:
\begin{displaymath}
 - \mathrm{div}(a(|x|)|\nabla w(x)|^{p-2} \nabla w(x)) + h(|x|,w(x),\langle\nabla w(x),\frac{x}{|x|}\rangle)=\phi(w(x)), 
\end{displaymath} 
where $x$ belongs to the ball in ${\R}^n$. Using the argument based on Opial-type inequalities, we investigate qualitative properties of their radial solutions, like e.g. maximum principles, monotonicity, as well as nonexistence of the nontrivial solutions.
\end{abstract}

\emph{Keywords}:
35B50, %maximum principles
%35B65 smoothness and regularity
35J92, %Quasilinear elliptic equations with $p$-Laplacian
%35J70 Degenerate elliptic equations
26D10. %Inequalities involving derivatives and differential and integral operators
\section*{Introduction}
We are interested in qualitative properties of  radial solutions of the following PDE: 
\begin{equation}\label{gendiv}
  - \mathrm{div}(a(|x|)|\nabla w(x)|^{p-2} \nabla w(x)) + H\Big(|x|,w(x),\langle\nabla w(x),\frac{x}{|x|}\rangle)\Big)=\phi(w(x)), 
\end{equation}
defined almost everywhere on a ball $B=B(0,R)\subseteq \mathbb{R}^n$, where $w \in W^{1,1}_{loc}(B)$, $|\nabla w|^{p-2}\nabla w \in L^1_{loc}(B,\R^n)$, $1<p<\infty$ and $a\in W^{1,1}_{loc}((0,R))$, $H:
[0,\infty)\times\R\times\R\rightarrow \R$, $\langle\cdot ,\cdot \rangle$ is the standard inner product in ${\R}^n$. 

The equation \eqref{gendiv} may be viewed as a generalization of a simple eigenvalue problem involving $p$-Laplace operator:
\begin{displaymath}
 -\mathrm{div}(|\nabla w(x)|^{p-2} \nabla w(x)) = \la |w(x)|^{p-2}w(x), \ \la >0,
\end{displaymath}
 or, more precisely, as the special variant of the more general eigenvalue problem involving the nonlinear gradient term:
\begin{displaymath}
-\Delta_{p,a(x)}w(x)+F(x,w(x),\nabla w(x))=\phi(w(x)), 
\end{displaymath}
involving the weighted  $p$-Laplacian: \[ \Delta_{p,\rho(x)}(w(x)):=\mathrm{div}(\rho(x)|\nabla w(x)|^{p-2}\nabla w(x)),\] with radial weight function $\rho(\cdot)=a(|\cdot|)$.

One of  our main results formulated in  Theorem
 \ref{twcole} gives the sufficient  conditions on the structure of 
 \eqref{gendiv} such that $|w(x)|$ has its supremum (possibly $\infty$) at zero. 
In another statement, Theorem \ref{twcopr}, we give the sufficient  conditions to deduce the nonexistence of the nontrivial radial $C^1$ solutions to (\ref{gendiv}). 

In our considerations, we proceed at first with the  equation in the nondivergent form:
\begin{multline}
- a(|x|)\mathrm{div}(|\nabla w(x)|^{p-2} \nabla w(x))
 + h(|x|,w(x),\langle\nabla w(x),\frac{x}{|x|}\rangle)=\phi(w(x)),\\
x\in B(0,R)\subset\R^n,\label{gen}
\end{multline}
where $h: [0,R)\times\R\times\R\rightarrow\R$.
It may be considered equivalent to \eqref{gendiv} when we require that
\begin{eqnarray*}
h(|x|,w(x),\langle \nabla w(x),\frac{x}{|x|}\rangle )~~~~~~~~~~~~~~~~~~~~~~~~~~~~~~~~~~~~~~~~~~~~~~~~~~\nonumber\\
= 
H( |x|,w(x),\langle \nabla w(x),\frac{x}{|x|}\rangle ) 
- a^{'}(|x|) |\nabla w(x)|^{p-2}\langle \nabla w(x),\frac{x}{|x|}\rangle , 
\end{eqnarray*}
that is when $H$ and $h$ are linked by
$$
h(s,p,q)=H(s,p,q)- a^{'}(s)|q|^{p-2}q.$$
Note that in the radial case $w(x)=u(|x|)$ for some scalar function $u$, we have $\langle \nabla w(x),\frac{x}{|x|}\rangle
= \langle u^{'}(|x|)\frac{x}{|x|} ,\frac{x}{|x|}\rangle =u^{'}(|x|)$,
$|\nabla w(x)|=|u^{'}(x)|= |\langle \nabla w(x),\frac{x}{|x|}\rangle | $. %Therefore when $H$ and $h$ are linked by
%$$
%h(s,p,q)=H(s,p,q)- a^{'}(s)|q|^{p-2}q,
%$$
%then  equation \eqref{gen}
% is equivalent to (\ref{gendiv}). 
Because of the radiality assumptions, \eqref{gen} reduces to the following ODE:
\begin{equation}
 a(\tau)(\Phi_p (u'(\tau)))' + (n -1)\frac{a(\tau)}{\t}\Phi_p (u'(\tau)) - h(\tau, u(\tau), u'(\tau)) + \phi(u(\tau)) = 0, \label{op}
\end{equation}
satisfied for a.e. $\tau \in B(0, R)$, where we use general notation, the same for every $k\in\mathbf{N}$: 
\begin{equation}\label{phi}
\Phi_p(\lambda) = |\lambda|^{p-2} \lambda\ \mathrm{for}\  \la\in{\R}^k\setminus\{ 0\}\ {\rm and}\  \Phi_p(0)=0. 
\end{equation}
Equation \eqref{op} is the starting point in our analysis. In particular, in Theorem \ref{prawa}, we contribute to the nonexistence results, as well as we provide the appriori estimates for the solutions.

The methods we use are further development of techniques from \cite{AdKa2009,AdKa2010,KaSt}, which were originated by Szeg\"o~\cite{Sze1} in the study of orthogonal polynomials. In particular in  \cite{KaSt}, the authors dealt with a linear variant of the equation (\ref{gendiv}) ($p=2$) and investigated some special functions like  Legendre, Jacobi polynomials,  Laguerre polynomials, or hypergeometric functions. The two subsequent papers \cite{AdKa2009,AdKa2010} focused on the application of that method to $p$-harmonic problems.  The authors have shown that, under some assumptions, the local maxima of the modulus of any radial solution form monotone sequence, which is a variant of the maximum principle. In \cite{AdKa2009} the authors deal with $h\equiv 0$, while in   \cite{AdKa2010}, in some  results  
it is assumed  that for a.e. $\t\in(0,R)\ \mathrm{and\ every}\ \la_0,\la_1\in\R$ the function $\h$ satisfies the following pointwise estimate (see Theorem 2.1):
\begin{equation}\label{hp}
h(\t,\la_0,\la_1)\la_1\leq\delta_a(\t)|\lambda_1|^p\ 
%\mathrm{or}\ h(\t,\la_0,\la_1)\la_1\geq\delta_a(\t)|\lambda_1|^p, 
\end{equation}
where
\begin{equation}
\delta_a (\t):=(n - 1)\frac{a(\t)}{\t}- \left(1-\frac{1}{p}\right)a'(\t) \ge 0\ a.e.. %\footnote{sprawdzic}. \label{delta}
\end{equation}
%It was not visible then how to study monotonicity property for radial solutions to (\ref{gen}) when the function $h(\cdot,\cdot,\cdot)$ satisfies some other estimate.
We contribute by proving similar type results when $\h$ satisfies different pointwise estimates:
\begin{equation*}
h(\t, \la_0 , \la_1)\la_1\leq q(\t)|\la_0|^l|\la_1|^{p-l}\ \mathrm{for}\ \mathrm{all}\ \la_0, \la_1 \in \mathbb{R}, \ \hbox{where}\  0<l<p, l\in\mathbf{R}, \label{innoh} 
\end{equation*}
  involving some nonnegative measurable function $q(\cdot)$ defined on $(0, R)$. It is related to $\delta_a(\cdot)$ via certain integral inequality (conditions from $\mathcal{A}5$ in Section \ref{ass-and-eq}). 
To prove our main results, instead of  pointwise inequality in \eqref{hp}, we use the Opial - type 
inequality  due to  Beesack and Das~\cite{BeeDas} (Theorem \ref{thmBeDa}), to deduce that  we now have its weaker, integral variant:
\begin{equation}
\int_a^b h(\t,u(\t),u'(\t))u'(\t)d\t\leq\int_a^b\delta_a(\t)|u'(\t)|^pd\t , \label{dem} 
\end{equation}
which appears sufficient for our analysis. In particular, Opial - type inequality,  
due to Bessak and Das,
 serves as a tool in the study of monotonicity properties for radial solutions to PDEs. To our best knowledge, such an application has not been noticed so far. 
We believe that the presented method, as well as Opial-type inequalities, can be further developed and applied  to the study of monotonicity properties of solutions to PDEs in the more general setting.

\smallskip
\noindent
Singular boundary value problems involving $p$-Laplacian arise for example, in fluid dynamics (\cite{CaNa}, \cite{Diaz}, \cite{Die}, Chapter~2 in \cite{Dra1}, \cite{WuZhLi}); glaciology (\cite{ArDiTe}), stellar dynamics (\cite{KaYaYo}); in the theory of electrostatic fields (\cite{FoOrPi}); in quantum physics (\cite{BDFP});  in the nonlinear elasticity theory (\cite{DOIw}).  
For some related topics, dealing with existence/nonexistence problem for singular or nonsingular PDE's we refer e.g. to \cite{bfg,cm,DraKuNi,filipucci1,pohmi_99_b} and to their references.

\section{Preliminaries} 
{\bf Notation.}  We will be dealing with the following spaces:
\begin{align*}\label{lokalne}
L^p_{loc}([0, R))&:=\bigcap_{0<r<R}L^p((0, r)),\ \ 
W^{n,p}_{loc}([0, R)):=\bigcap_{0<r<R}W^{n,p}((0, r)), 
\end{align*}
 as well as with their obvious analogues: $L^p_{loc}((0,R])$, $W^{n,p}_{loc}((0,R])$, where $L^p((0,r))$ and $W^{n,p}((0,r))$ are the usual $L^p$ and Sobolev spaces defined on the interval.
 We will also consider weighted variants of such spaces, like for example $L^{p}_{loc}((0,R], \rho dx)$,
 where $\rho$ is the weight.
 %, as well as weghted Beppo Levi type sets:
% $$
% \mathcal{L}^{1,p}((a,b)):=\{ u\in L^1_{loc}((a,b)): u^{'}\in 
 %L^p((a,b),\rho dx)  \}, 
 %$$ 
 %where $u^{'}$ is distrbutional derivative, and their local variants, %e.g. $\mathcal{L}_{loc}^{1,p}([a,b))$.  

Moreover, we deal with $\Phi_p(\cdot)$ as in~(\ref{phi}). 
%\begin{equation}
%\Phi_p(\lambda) = |\lambda|^{p-2} \lambda\ \mathrm{for}\ \la\neq 0, \Phi_p(0)=0, \la\in\R. %\label{phi}
%\end{equation}

\smallskip
\noindent
{\bf The variant of Opial-type inequality due to Beesack and Das.} We will use the following variant of Opial-type inequality due to Beesack and Das \cite{BeeDas} (see also e.g. \cite{Lee,PachD87,sin91,Shum} for some later contributions).
\begin{tw}
\label{thmBeDa}[Beesack and Das, 1968]\label{besacdas}. Let $l, m$ be real numbers such that $l,m>0$ and $l+m>1$, $-\infty< a < b< \infty$ and  
$p(\cdot),q(\cdot)$
be non-negative, measurable functions defined on  $(a,b)$ such that 
\begin{equation}\label{pep}
\int_a^b (p(t))^{-\frac{1}{l+m-1}} d t < \infty.
\end{equation}
 Assume further that the quantity $K(y)=K(a,y,l,m,q(\cdot),p(\cdot))$,   defined  for $a\leq y \leq b$ by
\begin{displaymath}
K(y)\!:=\!\left(\frac{m}{l+m}\right)\!^{\frac{m}{l+m}}\!\left[\int_a^y\! q(t)^{\frac{l+m}{l}} p(t)^{-\frac{m}{l}}\!\left(\int_a^t\! p(s)^{-\frac{1}{l+m-1}}ds\!\right)\!^{l+m-1}dt \right]^\frac{l}{l+m}, 
\end{displaymath}
  is finite. If $u(\cdot)$ is absolutely continuous on $[a,y]$ and either $u(a)=0$ or $u(y)=0$, then:
\begin{equation}\label{kropka}
\int_a^y q(t)|u(t)|^l|u'(t)|^{m} dt \leq K_1(y)\int_a^y p(t) |u'(t)|^{l+m} dt. 
\end{equation}
Moreover, the equality holds for all $y\in [0,b]$
if and only if either $u\equiv 0$  or there exist some constants $k_1\geq 0, k_2\in\mathbb{R}$ such that for $v(y):=\int_a^y p^{-[\frac{1}{l+m-1}]}dt$ we have: 
\begin{displaymath}
u(y)=k_2 v(y)\ \mathrm{and}\ q(y)=k_1 p^\frac{m-1}{l+m-1}\left(v(y)\right)^\frac{l(1-m)}{m}. 
\end{displaymath}
\end{tw}
%\footnote{sprawdzic}
%\begin{re}

\vspace{2em}

\noindent
Theorem \ref{thmBeDa} is a weighted variant of the classical Opial inequality~\cite{Opial}: $$\int_0^b |u(t)u'(t)|dt \leq \frac{b}{4}\int_0^b |u'(t)|^2 dt,$$ which holds for $u\in C^1(0,b)\cap C[0,b]$ such that $u(0)=u(b)=0$ and $u(t)>0$
for every $t\in(0,b)$. However, in the above statement we deal with different boundary conditions for $u$ and positivity of $u$ is not required.
%\end{re}
\section{Results for ODE's} \label{dowod}
\subsection{Assumptions and the associated equations}\label{ass-and-eq}
We will consider the following set of assumptions:

\smallskip
\noindent
{\bf General assumptions.}
\begin{enumerate}[label=$\mathcal{A}$\arabic*]
\item (Assumptions on the involved numbers) $p\in(1,\infty), l\in (0,p), n\in [1,\infty)$, $R\! \in\! (0, \infty]$. If $R\! =\! \infty$, then $(0,R)$ denotes the whole $\mathbb{R_+}$. 
\item (Assumptions about $\phi(\cdot)$) $\phi:\mathbf{R}\rightarrow \mathbf{R}$ is continuous, odd function  and we consider the following additional assumptions:
\begin{enumerate}
\item 
$\t \phi(\t) > 0$ for\ almost every  $\t;$
\item $\t \phi(\t) < 0$ for almost every $\t$.
\end{enumerate}
\item (Assumptions about $a(\cdot)$) $a(\cdot)\in W^{1,1}_{loc}((0, R))$, $a> 0$ a.e. (ellipticity condition).
 The functions $\delta_a(\cdot)$ and $d_a(\cdot)$ are defined a.e. on $(0,R)$ by:
\begin{eqnarray*}
\delta_a (\t)&:=&(n - 1)\frac{a(\t)}{\t}- \left(1-\frac{1}{p}\right)a'(\t);\\
d_a (\t)&:=&(n - 1)\frac{a(\t)}{\t} + \frac{1}{p}a'(\t).
\end{eqnarray*}
Moreover,  for an interval $X\in \{ (0,R), [0,R), (0,R] \}$ and $v(\cdot)\in \{ \delta_a(\cdot ), d_a(\cdot)\}$, which   is  positive  a.e., ,  we consider the additional assumptions:
\begin{equation}\label{wyjazd}
 a(\cdot)\in W^{1,1}_{loc}(X),\ 
 \ v(\cdot),v(\cdot)^{-1/(p-1)}\in L^1_{loc}(X),\ {\rm where}
\end{equation}
\begin{enumerate}
% \item[(a)]  $X=(0,R)$ and $v:=\delta_a$;
%\label{a}
 \item[($a_l$)]  $X=[0,R)$ and $v:=\delta_a$;
 \item[($a_r$)]  $X=(0,R]$ and $v:=\delta_a$;
\label{al}
\item[(a)]  $X=(0,R)$ and $v:=\delta_a$;
\item[($b_l$)]  $X=[0,R)$ and $v:=d_a$;
\item[($b_r$)]  $X=(0,R]$ and $v:=d_a$;
\item[(b)]  $X=(0,R)$ and $v:=d_a$.
\end{enumerate}
Note that in particular:  $a\in C([0,R))$ in case of $(a_l),(b_l)$;\\ 
 $a\in C((0,R])$ in case of $(a_r),(b_r)$;\\   $a\in C((0,R))$ in case of $(a),(b)$.
\item (Assumptions about $\h$) \label{anh} $h = h(\tau, \lambda_0 ,\lambda_1 ): (0, R)\times \mathbb{R}^2 \rightarrow \mathbb{R}$ is a Carath\'{e}odory function, i.e. is measurable with respect to $\tau\in (0,R)$ and continuous with respect to the remaining variables. Additionally, 
\begin{eqnarray}
%|h(\t, \la_0 , \la_1)|&\leq& \theta q(\t)|\la_0|^l|\la_1|^{p-l-1} + (1-\theta)v(\t)|\la_1|^{p-1},\ \hbox{\rm so that}\nonumber\\
h(\t, \la_0 , \la_1)\la_1&\leq& \theta q(\t)|\la_0|^l|\la_1|^{p-l} + (1-\theta)v(\t)|\la_1|^p, \label{a4}
\end{eqnarray}
for all $\la_0, \la_1 \in \mathbb{R}$, where $\theta \in (0,1], l\in (0,p)$,  $v(\cdot), q(\cdot): (0,R)\rightarrow [0,\infty)$ are nonnegative a.e. and belong to $L^1_{loc}(X)$ and 
%in case of $\alpha\neq 1$
 we deal with  $X$ and  $v$
as in $(a_l),(a_r),(a),(b_l),(b_r),(b)$ from $\mathcal{A}3$.  
In case of $\theta =1$ the assumptions on $v$ can be omited.

\item (Relation between $q(\cdot)$ and $a(\cdot)$) 
The following quantity is  defined in terms of the nonnegative a.e. measurable function 
$q: (0,R)\rightarrow [0,\infty)$, strictly positive a.e. function $v: (0,R)\rightarrow (0,\infty)$ ($v$, and parameters   $0\le s,r<R$, 
$1<p<\infty$, $l\in (0,p)$:
\begin{eqnarray*}
K(s, r,q,v):=\!\!\left(\frac{p-l}{p}\right)\!^{\frac{p-l}{p}}\!\left[\int_s^r q(t)^{\frac{p}{l}} v (t)^{-
\frac{p-l}{l}}\!\left(\int_s^t v(s)^{-\frac{1}{p-1}}ds\right)\!^{p-1} dt \right]^\frac{l}{p}\! .
\end{eqnarray*}
Moreover, we consider the following set of conditions linking $q(\cdot)$ and $a(\cdot)$, where $a(\cdot)$, $\delta_a(\cdot)$,  $d_a(\cdot)$ are as  in $\mathcal{A}3$:
\begin{enumerate}
 \item[($a_l$)] 
 $K(0, r,q,\delta_a)
< \infty\ \hbox{\rm for every}\ 0<r<R;$
 \item[($a_r$)] 
$K(0, R,q,\delta_a)
\le 1$;
\item[(a)] 
 $K(s, r,q,\delta_a)
< \infty\ \hbox{\rm for every}\ 0<s<r<R;$
\item[($b_l$)] 
$K(0, r,q, d_a)< \infty\ \hbox{\rm for every}\ 0<r<R;$
\item[($b_r$)] 
$K(0, R,q, d_a)\le 1$;  
\item[(b)] 
$K(s, r,q, d_a)< \infty\ \hbox{\rm for every}\ 0<s<r<R.$
\end{enumerate}
Note that in particular:\\
$\delta_a^{-1/(p-1)}\in L^1_{loc}([0,R))$ in case of $(a_l),(a_r)$;\\
$\delta_a^{-1/(p-1)}\in L^1_{loc}((0,R))$ in case of $(a)$;\\
$d_a^{-1/(p-1)}\in L^1_{loc}([0,R))$ in case of $(b_l),(b_r)$;\\
$d_a^{-1/(p-1)}\in L^1_{loc}((0,R))$ in case of $(b)$.

\end{enumerate}
%\begin{re}

\noindent
{\bf Remarks about the assumptions.}\\

\begin{re}\label{bp}
 Let us discuss the condition $\mathcal{A}3$.\\
 1) When the condition $\mathcal{A}3(a_l)$ holds we have: 
$$\delta_a(t)=\frac{n-1}{t}a(t)-\left(1-\frac{1}{p}\right)a'(t)> 0\ \ {\rm a.e.} \Leftrightarrow\ a'(t)\leq \frac{n-1}{t(1-\frac{1}{p})} a(t) \ \ {\rm a.e.} ,$$%TODO kiedy a jest 0
Gronwall's lemma yields:
$$a(t)< a(t_0) \left(\frac{t}{t_0}\right)^\frac{n-1}{(1-\frac{1}{p})}\ \mathrm{for\  all}\ t>t_0>0.$$ In particular the function $a(t)/t^\alpha$ is strictly dereasing on $(0,R)$ for  $\alpha=\tfrac{n-1}{(1-\frac{1}{p})}$ and $a(t)/ t^\alpha\ge C>0$  near zero.
\\
2) Similar estimates applied to $d_a$ in $\mathcal{A}3(b)$,
give 
$$a(t)< a(t_0) \left(\frac{t}{t_0}\right)^{-(n-1)p}\ \mathrm{for\  all}\ t>t_0>0$$
and consequently 
$a(t)t^{(n-1)p}$ is strictly decreasing and  $a(t)t^{(n-1)p}\ge C>0$ near zero.\\ 
3) The condition $\int_a^r  \rho(\tau)^{-1/(p-1)}d\tau <\infty$ (the local variant of \eqref{pep}) with the positive a.e., measurable function $ \rho(\cdot)$ and $a<r< b$, implies 
that\\ $L^p((a,r),  \rho(\cdot)dx)\subseteq L^1((a,r))$. This easily follows from H\"older's inequality and the observation comes from \cite{kuf-opic}:
\begin{eqnarray*}
\int_a^r |h(\tau)|\, d\tau &=& \int_a^r  (|h(\tau)| \rho(\tau)^{\frac{1}{p}})  \rho(\tau)^{-\frac{1}{p}} \, d\tau \\
&\le& 
\left( \int_a^r  |h(\tau)|^p  \rho(\tau)\, d\tau  \right)^{\frac{1}{p}} \left(  \int_a^r    \rho(\tau)^{-\frac{1}{p-1}}\, d\tau   \right)^{1-\frac{1}{p}} <\infty .
\end{eqnarray*}
We can apply this observation to $ \rho=\delta_a$ or $ \rho=d_a$, respectively
because of conditions from $\mathcal{A}3$. 
The condition  $\rho^{-1/(p-1)}\in L^1_{loc}([a,b))$   is stronger than the $B_p$  condition due to Kufner and Opic (\cite{kuf-opic}), where one assumes that $\rho^{-1/(p-1)}\in L^1_{loc}((a,b))$. 
\end{re}

\begin{re}\rm
Observe that the function $\phi(\cdot)$ in $\mathcal{A}2(a)$ changes its sign at $0$, as it is even.  
\end{re}

\begin{re}\rm
The number $n\ge 1$ in $\mathcal{A}1$ serves as an arbitrary real parameter. It will be interpreted as the dimension when we consider the multidimensional case.
\end{re}

\begin{re}\rm
The estimate \eqref{a4} with  $\theta =0$  was considered in \cite{AdKa2010}.
\end{re}

\smallskip
\noindent
{\bf The associated ODEs.}
We consider the following ODE's satisfied for a.e. $\t\in(0,R)$, having  nondivergent and divergent forms, respectively:
\begin{equation}
 a(\tau)(\Phi_p (u'(\tau)))' + (n -1)\frac{a(\tau)}{\t}\Phi_p (u'(\tau))) - h(\tau, u(\tau), u'(\tau)) + \phi(u(\tau)) = 0 \label{o}
\end{equation}
and:
\begin{equation}
(a(\tau)\Phi_p (u'(\tau)))' + (n -1)\frac{a(\tau)}{\t}\Phi_p (u'(\tau))) - h(\tau, u(\tau), u'(\tau)) + \phi(u(\tau)) = 0, \label{do}
\end{equation}
The ODE (\ref{o}) is associated to the PDE (\ref{gen}) restricted to its radial solutions, 
%, which was staded in nondivergent form, 
while the equation (\ref{do}) is associated to the PDE (\ref{gendiv}) restricted to its radial solutions. 
The assumption $a>0$ a.e. interprets that the equation (\ref{o}) is elliptic. It becomes degenerate when $a(\cdot)$ achieves 0. 
%, which is statet in nondivergent form. 

In the preceding sections we will discuss the appriori estimates, nonexistence/triviality and monotonicity of their solutions.

\subsection{Nonexistence and triviality of  solutions}
We start with the presentation of nonexistence/triviality results. 

\subsubsection{Formulation of results}
We will deal with the following sets of conditions:
\begin{eqnarray*}
\mathcal{N}_{nd}&:=&
\{  \mathcal{A}1, \mathcal{A}2(a),\mathcal{A}3(a_l), \mathcal{A}4(a_l),  \mathcal{A}5(a_l)\},\\
\mathcal{N}_{d}&:=& \{\mathcal{A}1,\mathcal{A}2(a),\mathcal{A}3(b_l), \mathcal{A}4(b_l),\mathcal{A}5(b_l)\}. 
\end{eqnarray*}

Our first result is applicable to the equations  (\ref{o}) and (\ref{do}).  It is complementary to the result from \cite{AdKa2010}, where the authors have considered the solutions to (\ref{o}) and assumed the estimate
\begin{equation}\label{pop}
h(\tau,\lambda_0,\lambda_1)\lambda_1\le \delta_a(\tau)|\lambda_1|^p.
\end{equation}
Now we deal with the estimate \eqref{a4}, where $v\in\{ \delta_a,d_a\}$. Note that the case of $\theta= 0$, $v=\delta_a$ in \eqref{a4}, is precisely \eqref{pop}.  

In the statement below we obtain the nonexistence and triviality results for solutions to (\ref{o}) and (\ref{do}), as well as their appriori estimates in $L^\infty$, obtained  in terms of the boundary data.
The appriori estimates were not considered in  \cite{AdKa2010}.

\begin{tw}[Estimates and triviality]\label{prawa}~\\
Let   $u(\cdot)\in W_{loc}^{1,1}((0,R))$ be such that
%is a solution to the ODE (\ref{o}) such that
\begin{enumerate}[label=(\alph*)]
\item (Regularity assumption) 
%$u(\cdot) \in W^{1,1}_{loc} ((0,R))$, 
$\Phi_p(u'(\cdot)) \in W_{loc}^{1,1}([0,R))$, in particular $u,u^{'}\in C([0,R))$;
\item (Boundary condition) $u(0)=0$. 
%and $\lim_{\epsilon\to 0} a(\epsilon)|u'(\epsilon)|^p=0$;
\end{enumerate}
Moreover, suppose that  one of the following assumptions are satisfied:
\begin{itemize}
\item[($\mathcal{ND}$):]
 $u$ is a solution to the ODE \eqref{o}, $\mathcal{N}_{nd}$ holds, $K:= K(0,R,q,\delta_a)$;
\item[($\mathcal{D}$):] 
 $u$ is a solution to the ODE \eqref{do}, $\mathcal{N}_{d}$ holds, $K:= K(0,R,q,d_a)$.
\end{itemize} 
Then we have:
\begin{description}
\item[i)] When  $K\le 1$ 
and $\Phi(\t):=\int_0^\t \phi(s)ds$, then   
\begin{equation*}\label{gorac}
{\rm sup}_{r\in (0,R)}|u(r)|\le \Phi^{-1}\left( \left(1-\frac{1}{p}\right)a(0)|u'(0)|^p \right). 
\end{equation*}
\item[ii)]
When  $ a(0)=0$ or $u'(0)=0$, then the only solution can be $u\equiv 0$. In that situation and when $h(\cdot ,0,0)\not\equiv 0$ a.e.,
there are no such solutions.
\end{description}
\end{tw}

\subsubsection{Proof of Theorem \ref{prawa}}

{\bf Proof of Theorem \ref{prawa} under the assumption ($\mathcal{ND}$).}
 We use similar arguments as in the proof of Theorem 2.1 in  \cite{AdKa2010}. Let 
 $$ K(r) := K(0,r,l,p-l,q(\cdot),\delta_a(\cdot)) $$ 
where $K(\cdot)$ is as in $\mathcal{A}5(a_l)$. We have  $K(r)<\infty$ 
for every $0<r<R$.

\smallskip
\noindent 
{\sc Proof of part i):}  We assume  that $K(R)\leq 1$.\\
 Define:
\begin{align}
A(\t_1,\t_2)&:=\Phi(|u(\t_1)|)-\Phi(|u(\t_2)|).\nonumber
\end{align} 
We notice that due to condition $\mathcal{A}2(a)$ the function $\Phi\colon[0,R)\rightarrow[0,\infty)$ is strictly increasing and $\Phi(0)=0$. 
Since $\Phi$ is locally Lipschitz, it follows that $\Phi\circ|u|\in W_{loc}^{1,1}([0,R))$ and 
for every $r$ such that $ 0<r<R$ we have:
\begin{align*}
A(r,0)&=\int_0^r\frac{d}{dt}\Phi(|u(t)|)dt=\int_0^r \Phi'(|u(t)|)\mathrm{sgn}(u(t))u'(t)dt =\\
&=\int_0^r \phi(u(t))u'(t)dt.
\end{align*}
Note that $\phi(u)u^{'}$ is integrable over $(0,r)$.
Multiplying \eqref{o} by $u'$, then using  the pointwise estimate \eqref{a4},  we obtain:
\begin{eqnarray}
\phi(u(\t))u'(\t )
&\le & -a(\t)(\Phi_p(u'(\t)))'u'(\t) - (n - 1)\frac{a(\t)}{\t}|u'|^p \nonumber\\
%+ h(\t, %u(\t), u'(\t))u'(\t)\\
%&\le& -a(\t)(\Phi_p(u'(\t)))'u'(\t) - (n - 1)\frac{a(\t)}{\t}|u'|^p \\
&+&
 \theta q(\t)|u(\t)|^l|u'(\t)|^{p-l} + (1-\theta)\delta_a (\t)|u'(\t)|^p
  \label{zup}
\end{eqnarray}
a.e. on $(0,R)$, where each involved summand on the right hand side above is integrable over $(0,r)$, because $u(\cdot),u^{'}(\cdot),a(\cdot)$ are bounded on $(0,r)$ (assumption (a) and $\mathcal{A}3(a_l)$), while 
$(\Phi_p(u^{'}))^{'}$,  $q(\cdot)$, $\delta_a(\cdot)$, $a(\t)/\t$ are integrable over $(0,r)$ (by (a), $\mathcal{A}3(a_l)$, $\mathcal{A}4(a_l)$, and because $a(\t)/\t$ is represented in terms of $a^{'}$ and $\delta_a$). 

Integrating the above equation over $(0,r)$, then applying the Opial-type inequality (Theorem~\ref{thmBeDa})
with parameters: $a=0,y=r,l, m=p-l$ and functions $q(\cdot),$ $p(\cdot)=\delta_a(\cdot)$, together with the assumptions $\mathcal{A}3(a_l)$, $\mathcal{A}5(a_l)$, we get:
\begin{gather}
\int_0^r \theta q(\t)|u(\t)|^l|u'(\t)|^{p-l}d\t + \int_0^r (1-\theta)\delta_a(\t)|u'(\t)|^p d\t\leq 
\nonumber\\ \leq \theta K(R)\int_0^r\delta_a(\t)|u'(\t)|^p d\t + (1-\theta) \int_0^r \delta_a(\t)|u'(\t)|^p d\t \nonumber\\
%=  \bar{K}(r)\int_\epsilon^r\delta_a(\t)|u'(\t)|^p d\t 
\le \int_0^r\delta_a(\t)|u'(\t)|^p d\t  .\nonumber
\end{gather}
%where $\bar{K}(r)= \theta K(r) + (1-\theta)\le 1$.
\noindent
Altogether give:
\begin{gather*}
\int_0^r \phi(u(\t))u'(\t )d\t \leq  -\int_0^r a(\t)(\Phi_p(u'(\t)))'u'(\t)d\t - \int_0^r(n - 1)\frac{a(\t)}{\t}|u'(\tau)|^p d\t\\+\int_0^r \delta_a(\t) |u'(\t)|^p d\t.
\end{gather*}
Using the definition of $\delta_a(\cdot)$ from $\mathcal{A}3$ and the fact that each integrand on the right hand side above is integrable,  we get:
\begin{equation}
\begin{split}
\int_0^r &\phi(u(\t))u'(\t )d\t \leq \\ \leq \int_0^r & \left( -  a(\t)(\Phi_p(u'(\t)))'u'(\t) - \left(1-\frac{1}{p}\right)a'(\t)|u'(\t)|^p \right) d\t . \label{blko}
\end{split}
\end{equation}
We will estimate the right hand side of the inequality with the help of
%$\Psi(\t,\la_1):=-\left(1-\frac{1}{p}\right)a(\t)|\la_1|^p$, so that
\begin{equation}\label{prezentacja}
\Psi(\t,\la_1):=-\left(1-\frac{1}{p}\right)a(\t)|\la_1|^p,\ \hbox{\rm  i.e.}\ 
\Psi(\t, u'(\t))=-\left(1-\frac{1}{p}\right)a(\t)|u'(\t)|^p. 
\end{equation}
We note that because 
%$u^{'}\in L^1_{loc}((0,R))$ and 
$\Phi_p(u^{'})\in W^{1,1}_{loc}([0,R))$, it is absolutely continuous on $[0,R)$, while $v\mapsto |v|^{p/(p-1)}$ is locally Lipschitz. Therefore also   $|u^{'}|^p=|\Phi_p(u^{'})|^{\frac{p}{p-1}}$ is absolutely continuous on $[0,R)$ and we can use  differentiation formula for compositions:
\begin{eqnarray}
(|u'|^p)'=\frac{p}{p-1}u'\cdot(\Phi_p(u'))^{'}\ {\rm a.e.,\  consequently}\nonumber\\
\frac{d}{d\tau} \Psi (\tau, u^{'}(\tau))= - a(\t)(\Phi_p(u'(\t)))'u'(\t) 
- (1-\frac{1}{p})a^{'}(\t)|u'(\tau)|^p. \label{dodatk}
\end{eqnarray}
Hence and from \eqref{blko}, for any $r\in (0,R)$: 
%the right hand side in (\ref{blko}) is equal to  
\begin{eqnarray*}\label{pstryk}
\Phi (|u(r)|)= A(r,0)\le \int_0^r\frac{d}{d\t} \Psi(\t,u'(\t))d\t =\Psi(r, u^{'}(r))-\Psi(0, u^{'}(0))\nonumber\\
\le  \left(1-\frac{1}{p}\right)a(0)|u'(0)|^p . 
\end{eqnarray*} 
 This gives the estimate in part i).

\smallskip
\noindent
{\sc Proof of part ii):} We have no restrictions on $K(\cdot)$ except its  finiteness on $(0,R)$.\\
We can thus assume that $1\le K(r)<\infty$, as otherwise the conclusion follows from already proven part i).
Consider
$$
I:= \{ 0\} \cup\{ r\in (0,R): u\equiv 0\ {\rm on}\ [0,r)\} .
$$
Obviously  $I\neq \emptyset$. We easily verify that $I$ is connected and  closed, because of the continuity of $u$. To finish the proof of assertion ii) it suffices to show that $I$ is also open, as then $I=[0,R)$.
For this,  let $r_0\in I$. Because of the assumption $\mathcal{A}5(a_l)$, we will find $\rho\in (r_0,R)$ such that
$K(r_0,\rho,q(\cdot),\delta_a(\cdot),l,p-l)\le 1$. We can thus examine $A(r,r_0)$ with $r\in (r_0,\rho)$ with minor changes in the computations in the proof part i), 
where we change the integrals over $(0,r)$ by the ones over $(r_0,r)$
%, use $K(r_0,\rho,l,p-l,q(\cdot),\delta_a(\cdot))$ instead of previous $K(R)$ 
and remember that   $u(r_0)=u^{'}(r_0)=0$.
They show that $u\equiv 0$ on $(r_0,\rho)$. This implies openness of $I$ and then the triviality/nonexistence assertions under ($\mathcal{ND}$) .

\subsubsection*{Proof of Theorem \ref{prawa}  under the assumption ($\mathcal{D}$)}
Equation \eqref{do} is equivalent to the nondivergent one:
\begin{eqnarray}\label{rownoniediw}
 a(\tau)(\Phi_p (u'(\tau)))' + \left\{ (n -1)\frac{a(\tau)}{\t} + a^{'}(\tau) \right\}\Phi_p (u'(\tau))~~~~~~~~~~~~~~~~~~~~~~~~ \\
 ~~~~~~~~~~~~~~~~~~~~~~~~~~~~~~~~- h(\tau, u(\tau), u'(\tau)) + \phi(u(\tau)) = 0. \nonumber
\end{eqnarray}
%Let us set
%$$
%\tilde{h}(\tau,\lambda_0,\lambda_1):= h (\tau,\lambda_0,\lambda_1)+ %b(\tau)\Phi_p(\lambda_1).
%$$
%Note that 
%$$
%\tilde{h}(\tau,\lambda_0,\lambda_1)\lambda_1\le q(\tau)|\lambda_0|^l|%\lambda_1|^{p-l}
%$$
We provide almost the same proof as that of Theorem \ref{prawa} with the following modifications in the proof of part i):
\begin{itemize}
\item now
we deal with 
$$
K(r):= K_1(0,r,l,p-l, q(\cdot), d_a(\cdot)),
$$
which is finite for all $0<r<R$ due to the assumption $\mathcal{A}5(b_l)$;
\item instead of \eqref{zup}  we have:
\begin{eqnarray*}
 \phi(u(\t))u'(\t )d\t \le\\
  - a(\t)(\Phi_p(u'(\t)))'u'(\t)d\t  
 +\left\{ (n - 1)\frac{a(\t)}{\t}
 + a^{'}(\tau)\right\}|u'|^p d\t\\  
 + \theta q(\t) |u(\t)|^l|u^{'}(\t)|^{p-l} + (1-\theta )d_a(\t)|u^{'}(\t)|^p.
\nonumber
\end{eqnarray*}
\end{itemize}
By the modified computations we arrive at \eqref{blko} and the remaining arguments are the same as in the proof of the case under ($\mathcal{ND}$).
%That ends the proof of Theorem \ref{prawa}. 
\hfill$\Box$
\subsubsection{Remark about the support of the solution}
We will consider the following sets of conditions:
\begin{eqnarray*}
\mathcal{S}_{nd}&:=&
\{  \mathcal{A}1, \mathcal{A}2(a),\mathcal{A}3(a), \mathcal{A}4(a),  \mathcal{A}5(a)\},\\
\mathcal{S}_{d}&:=&
\{  \mathcal{A}1, \mathcal{A}2(a),\mathcal{A}3(b), \mathcal{A}4(b),  \mathcal{A}5(b)\},
\end{eqnarray*}

\noindent
Following  the proof of Theorems \ref{prawa}, where the only difference is that we assume $u(s)=u^{'}(s)=0$ for some $s
\in (0,R)$, and we validate the integrals over $(s,r)$ for $r\in (s,R)$, we obtain the following statement. Its similar proof is left to the reader.

\begin{tw}[Support of the solution]\label{support}~\\
Let   $u,\Phi_p(u'(\cdot)) \in W_{loc}^{1,1}((0,R))$ and either i) or ii) is satisfied, where:
\begin{description}
\item[($\mathcal{SND}$)] $\mathcal{S}_{nd}$ holds,  $u(\cdot)$ is a solution to the ODE (\ref{o});
\item[($\mathcal{SD}$)] $\mathcal{S}_{d}$ holds,  $u(\cdot)$ is a solution to the ODE (\ref{do}).
\end{description}
Then if $u$ tauches zero at some its critical point $s\in (0,R)$, then  $u\equiv 0$  on $[s,R)$. 
\end{tw}

\vspace{1em}
\noindent
%The examples are postponed to the last sections.

\subsection{Left hand side maximum principles and monotonicity}
In this section we apply Opial-type inequality from Theorem \ref{besacdas} to function $u(\cdot)$, which is equal to zero at the right end of its domain. It allows  to deduce  that $u$ is of constant sign and monotone. In particular,  $|u|$ achieves its supremum at $0$ (possibly infinite). Precise formulation is given below.

\subsubsection{Formulation of results}
Consider  the following conditions:
\begin{eqnarray*}
\mathcal{M}_{nd} &:=& \{ \mathcal{A}1, \mathcal{A}2(b), \mathcal{A}3(a_r), \mathcal{A}4(a_r),\mathcal{A}5(a_r)\} ,
  \\ 
\mathcal{M}_{d} &:=& \{ \mathcal{A}1,\mathcal{A}2(b),\mathcal{A}3(b_r),\mathcal{A}4(b_r),\mathcal{A}5(b_r)\} . 
\end{eqnarray*}

We obtain the following maximum principle.

\begin{tw}[Constant sign and monotonicity]\label{lewa}
Let $u(\cdot):(0,R)\rightarrow \R$ be such that
\begin{enumerate}[label=(\alph*)]
\item (Regularity assumption)
$u\! \in\! W_{loc}^{1,1} ((0,R))$,  $\Phi_p(u')\! \in\! W_{loc}^{1,1} ((0,R])$, in particular $u,u^{'} \in C((0,R])$; 
\item (Boundary condition)   $u(R)=0$. 
\end{enumerate}
Moreover, assume that one of the conditions $(\mathcal{MND})$ or $(\mathcal{MD})$ holds where:
\begin{description}
\item[($\mathcal{MND}$):]  $u(\cdot)$ is a solution to \eqref{o} and $\mathcal{M}_{nd}$ holds;
\item[($\mathcal{MD}$):] $u(\cdot)$ is a solution to \eqref{do} and $\mathcal{M}_{d}$ holds.
\end{description}
Then $u(\cdot)$ is of constant sign and monotone, moreover:
\begin{equation*}
\sup_{x\in(0,R)}|u(x)|=\lim_{\epsilon\rightarrow 0}|u(\epsilon)|.
\end{equation*} 
If additionally $u(0)=0$ or $\limsup_ {\epsilon\rightarrow 0} a(\epsilon)|u'(\epsilon)|^p = 0$ then there are no such nontrivial solutions. 
\end{tw} 

%\footnote{Problem Dirichleta ma tylko zerowe rozwiazanie}

\subsubsection{Proof of Theorem \ref{lewa}}
We present separately the proof under the assumptuion $(\mathcal{MND})$ and $(\mathcal{MD})$. 

\medskip
\noindent
{\bf Proof under the assumption ($\mathcal{MND}$).}\\
We  denote: 
\begin{displaymath} 
\begin{split}
\tilde{\phi}(\t):=-\phi(\t),\ \ \ 
\tilde{\Phi}(\t):=\int_0^\t \tilde{\phi}(s)ds,\ \ 
\tilde{A}(\t_1,\t_2):=\tilde{\Phi}(|u(\t_1)|)-\tilde{\Phi}(|u(\t_2)|),\\
\tilde{\Psi}(\t,\lambda_1):= (1-\frac{1}{p})a(\t)|\lambda_1|^p= -\Psi (\t,\lambda_1) \ \ \hbox{\rm (see\ \eqref{prezentacja}).}~~~~~~~~~~~~~~~~~~~~~~ 
\end{split}
\end{displaymath}
Note that, due to our assumption $\mathcal{A}2(b)$, $\tilde{\Phi}$ is strictly increasing and  $\tilde{\Phi}(0)=0$. 
Let 
$$
K(r):= K_1(r,R,l,p-l, a(\cdot),\delta_a(\cdot)), \ \ r\in (0,R), 
$$
where let $K(\cdot )$ is as in $\mathcal{A}5$. The assumption $\mathcal{A}5(a_r)$ implies $K(r)\le 1$ for every $0<r\le R$. 
From now the proof follows by steps.

\smallskip
\noindent
% We start with the proof of the monotonicity of $u(\cdot)$. \\
{\sc STEP 1:} We prove that for any critical point $r$ of $u(\cdot)$ we have $u(r)=0$. It is enough to prove that $\tilde{A}(R,r) \geq 0$ for such $r$'s, because then 
$0\ge -\tilde{\Phi}(|u(r)|)= \tilde{\Phi}(|u(R)|)-\tilde{\Phi}(|u(r)|)= \tilde{A}(R,r)\ge 0$, consequently $u(r)=0$. 
For this, we compute that:
\begin{gather}
\tilde{A}(R,r) = \tilde{\Phi}(|u(R)|)-\tilde{\Phi}(|u(r)|)=\int_r^R
 \tilde{\phi}(u(\t))u'(\t)d\t, \ {\rm where}\nonumber\\
\tilde{\phi}(u(\t))u'(\t) \stackrel{\eqref{zup}}\ge 
  a(\t)(\Phi_p(u'(\t)))'u'(\t)  + (n - 1)\frac{a(\t)}{\t}|u'(\t)|^p  
 \nonumber\\ 
-\theta q(\t) |u(\t)|^l|u^{'}(\t)|^{p-l} - (1-\theta) \delta_a(\t)|u^{'}(\t)|^p.\label{zolodz}
\end{gather}
Moreover, all the involved terms in the last two lines above are integrable over $(r,R)$. Integrating them over $(r,R)$,
 applying Theorem \ref{besacdas} with $a=r, b=R$, recalling 
that $u(R)=0$ due to (b), and that  $K(r)\le 1$,  we deduce that for every $r\in (0,R)$ (not necessarily being the critical point of $u$):  
\begin{gather}
\tilde{A}(R,r)\geq \int_r^R \left( a(\t)(\Phi_p(u'(\t)))'u'(\t) +(n - 1)\frac{a(\t)}{\t}|u'(\t)|^p - \delta_a(\t) |u'(\t)|^p \right) d\t \nonumber\\\stackrel{\mathcal{A}3}{=}\int_r^R \left( a(\t)(\Phi_p(u'(\t)))'u'(\t) + \left(1-\frac{1}{p}\right) a'(\t)|u'(\t)|^p\right) d\t \nonumber\\
 \stackrel{\eqref{dodatk}}{=}\int_r^R \frac{d}{d\t}\tilde{\Psi}(\t,u'(\t)) d\t.\label{kwiecien}
\end{gather}
The nonnegativity of $a(\cdot)$ allows us to conclude, that for $r$'s --- critical points of $u(\cdot)$ -- there holds:
\begin{gather*}
\tilde{A}(R,r)\geq \left(1-\frac{1}{p}\right)a(R)|u'(R)|^p - \left(1-\frac{1}{p}\right)a(r)|u'(r)|^p  \\
= \left(1-\frac{1}{p}\right)a(R)|u'(R)|^p \geq 0.
\end{gather*}
This completes the proof of {\sc Step 1}. \\
{\sc STEP 2:} We have proven that for any critical point $r$ of $u(\cdot)$ we have $u(r)=0$. Thus, and because of Fermat's Theorem, $u$ cannot have local extrema (positive nor negative) inside $[0,R]$, and so $u$ is monotone (not necessarily strictly).  
This completes the proof of first part of the statement.\\ 
{\sc STEP 3:} We prove last assertion of the statement, which is trivial when $u(0)=0$. \\Assume that additionally $\lim_{\epsilon\rightarrow 0} a(\epsilon)|u'(\epsilon)|^p = 0$. Then:
\begin{gather*}
0\ge -\limsup_{\epsilon\rightarrow 0}\tilde{\Phi}(|u(\epsilon)|)= \tilde{\Phi}(|u(R)|)-\limsup_{\epsilon\rightarrow 0}\tilde{\Phi}(|u(\epsilon)|)= \liminf_{\epsilon\rightarrow 0}\tilde{A}(R,\epsilon)  \\
\stackrel{\eqref{kwiecien}}{\geq}\liminf_{\epsilon\rightarrow 0} \tilde{\Psi}(\t,u'(\t))|^R_\epsilon%\geq 
\geq \left(1-\frac{1}{p}\right)\left( a(R)|u'(R)|^p-\limsup_{\epsilon\rightarrow 0}a(\epsilon)|u'(\epsilon)|^p\right)\\=\left(1-\frac{1}{p}\right) a(R)|u'(R)|^p\geq 0,
\end{gather*}
which proves that $\lim_{\epsilon\to 0} u(\epsilon)=0$ and together with the monotonicity of $u$ and boundary condition (b) shows that $u\equiv 0$ on $[0,R]$. This completes the proof of  Theorem \ref{lewa} under ($\mathcal{MND}$).

\smallskip
\noindent
{\bf Proof under the assumption ($\mathcal{MD}$).}\\
We use the  non-divergent equivalent ODE \eqref{rownoniediw}  instead of 
\eqref{do}, and adopt the proof for $\mathcal{(MND)}$
 with the following modifications:
\begin{itemize}
\item in the proof of {\sc Step 1} %in case of $\mathcal{C}_3$ 
we deal with 
$$
K(r):= K_1(r,R,l,p-l, q(\cdot), d_a(\cdot)),
$$
which does not exceed $1$ for all $0<r<R$, due to  $\mathcal{A}5(b_r)$;
\item instead of \eqref{zolodz}  we have:
\begin{equation}
\begin{split}
 \tilde{\phi}(u(\t))u'(\t )\ge a(\t)(\Phi_p(u'(\t)))'u'(\t) + \left\{ (n - 1)\frac{a(\t)}{\t}+ a^{'}(\tau)\right\}|u'(\t)|^p 
 \\-  \theta q(\t) |u(\t)|^l|u^{'}(\t)|^{p-l} 
 - (1- \theta) d_a(\t) |u^{'}(\t)|^{p} .\nonumber
\end{split}
\end{equation}
\end{itemize}
Easy details are left to the reader.\hfill$\Box$

\section{Results for PDE's}
In this section we are interested in the multi-dimensional case. 

\subsection{The associated PDE's and auxilary fact}
We will deal now with radial solutions to PDE's:
\begin{eqnarray}\label{zabka}
- a(|x|)\mathrm{div}(|\nabla w(x)|^{p-2} \nabla w(x)) + h(|x|,w(x),\langle\nabla w(x),\frac{x}{|x|}\rangle)=\phi(w(x)),\\
- \mathrm{div}(a(|x|)|\nabla w(x)|^{p-2} \nabla w(x)) + h(|x|,w(x),\langle\nabla w(x),\frac{x}{|x|}\rangle)=\phi(w(x)),\label{zabkadiv}
\end{eqnarray}
in $\mathcal{D}^{'}(B)$,  where $u\in W^{1,1}_{loc} (B )$, $B=B(0,R)\subseteq \R^n$ is a ball, $n>1$, under certain assumptions, which will be discussed later.
 
\smallskip

The following lemma will be helpful to understand the interplay between regularity conditions related to the  multidimensional  case and that related to the one-dimensional case. It is obtained as a modification of 
Fact 2.1 from \cite{AdKa2010}. However, some arguments the proof of parts 4 and 5 do not follow so directly from previous ones. Therefore we present them  for reader's cenvenience. 

\begin{lem} \label{lemreg}
Let $n > 1, p > 1, w(x) = u(|x|)$ and \\
\mbox{$w\in W^{1,1}_{loc} (B \setminus \{0\})$}, where $B=B(0,R)\subseteq\mathbf{R}^n$,  $R\in\mathbf{R}$, $n\ge 2$. \\Then:
\begin{enumerate}
\item $u\in W^{1,1}_{loc}((0, R))$;
\item If $w\in W^{1,1}(B\setminus B(0,r))$ for any $0<r<R$, then $u\in W^{1,1}_{loc}((0, R])$, in particular  $u\in C((0,R])$;
\item If $\Phi_p(\nabla w)\in W^{1,1}_{loc}(B\setminus\{0\})$, then $\Phi_p(u')\in W^{1,1}_{loc}((0, R))$, in particular  $u, u^{'}\in C((0,R))$;
\item If $\Phi_p(\nabla w)\in W^{1,1}(B\setminus B(0,r) )$ for any $0<r<R$, then $\Phi_p(u')\in W^{1,1}_{loc}((0, R])$, in particular $u, u^{'}\in C((0,R])$;
\item If $\Phi_p(\nabla w)\in W^{1,1}_{loc}(B, |x|^{-(n-1)}dx)$, then $\Phi_p(u')\in W^{1,1}_{loc}([0, R))$, in prticular and $u, u^{'}\in C([0,R))$;
\end{enumerate}
\end{lem}

In the above notation we sometimes omit the fact, that the consider function is vector valuable, like e. g. $\Phi_p(\nabla w)$. However,
in some condiderations we need to mention it.

\bigskip
\noindent \textbf{Proof (of Parts 4,5)}: By parts 1 and 3 we have $u\in W^{1,1}_{loc}((0,R))$ and $\Phi_p(u^{'})\in W^{1,1}_{loc}((0,R))$. 
We start with the proof of part 4. \\
We have to prove that $\Phi_p(u^{'})\in W^{1,1}_{loc}((0,R])$, as the remaining  statement follows from Sobolev's Embedding Theorem. We note that
$$
\Phi_p(\nabla w (x)) =|\nabla w (x)|^{p-2}\nabla w (x) =|u^{'}(|x|)|^{p-2} u^{'}(|x|)\frac{x}{|x|} =\Phi_p(u^{'}(|x|)) \frac{x}{|x|}.
$$
We compute that for almost every $x\in B(0,R)$ and every $0<r<R$:
\begin{eqnarray}
L^1(B\setminus B(0,r), \mathbf{R}^n\times \mathbf{R}^n)&\ni& \nabla \Phi_p(\nabla w (x))\nonumber\\ &=& \left( \Phi_p(u^{'}(\tau))\right)^{'}|_{\tau = |x|} \frac{x}{|x|}\otimes \frac{x}{|x|} +  \Phi_p(u^{'}(|x|))\nabla \left( 
\frac{x}{|x|} \right)\nonumber\\
&=& (\Phi_p(u^{'}))^{'}|_{|x|}v(x) + \Phi_p(u^{'})|_{|x|}w(x),\ {\rm where},\nonumber\\
w(x)&:=& \nabla \left( 
\frac{x}{|x|} \right)\perp \frac{x}{|x|}\otimes \frac{x}{|x|}=: v(x)\ {\rm a.e.}.\label{szalik}
\end{eqnarray}
To explain \eqref{szalik} we note that for $x\neq 0$ we have $v_{ij}(x)=\frac{x_ix_j}{|x|^2}$, 
$w_{ij}= \frac{\partial}{\partial x_i}\left( \frac{ x_j}{|x|}\right) =
\frac{1}{|x|}\left(\delta_{ij} - \frac{x_ix_j}{|x|^2}   \right)$ and
$$
\sum_{i,j}v_{ij}w_{ij} =\frac{1}{|x|^3} \sum_{i}\sum_j\left(\delta_{ij} - \frac{x_ix_j}{|x|^2}   \right) x_ix_j  =\frac{1}{|x|^3} \sum_{i}\left( 
x_i^2- x_i^2
\right) =0.
$$
Moreover,
%\begin{eqnarray*}
%w(\cdot)&\in& L^\infty (B\setminus (B(0,r)),\mathbf{R}^n\times \mathbf{R}^n)\ \hbox{\rm for every $0<r<R$,}\\ 
$
| v(x)| = 1\ {\rm for}\ x\neq 0.
$
%\end{eqnarray*} 
Therefore
\begin{eqnarray*}
L^1(B\setminus B(0,r))\ni \nabla \Phi_p(\nabla w (x))\cdot v(x)= \left( \Phi_p(u^{'}(\cdot))\right)^{'}|_{|x|}.
\end{eqnarray*}
It follows that 
\begin{eqnarray*}
\int_r^R |\left( \Phi_p(u^{'})\right)^{'} (\tau)| d\tau &\le& \frac{\theta_{n-1}}{\theta_{n-1}r^{n-1}}
\int_r^R |\Phi_p(u^{'})^{'} (\tau)| \tau^{n-1}\, d\tau \\
&=& 
\frac{1}{\theta_{n-1}r^{n-1}}\int_{B\setminus B(0,r)}  |\left( \Phi_p(u^{'})\right)^{'}(|x|) |dx <\infty ,
\end{eqnarray*} 
where $\theta_{n-1}$ is the
$n-1$-dimensional  Lebegue's measure of the unit sphere in $\mathbf{R}^n$. This ends the proof of part 4 in the lemma.

To prove part 5, we modify our last inequality to the following:
\begin{eqnarray*}
\int_0^r |\left( \Phi_p(u^{'})\right)^{'} (\tau)| d\tau &=& 
 \int_0^r \left(  |\left( \Phi_p(u^{'})\right)^{'} (\tau)|\tau^{-(n-1)}\right)\tau^{(n-1)} \, d\tau\\
&=&
\frac{1}{\theta_{n-1}}\int_{B(0,r)}  |\left( \Phi_p(u^{'})\right)^{'}(|x|) ||x|^{-(n-1)}\, dx <\infty .
\end{eqnarray*} 
\hfill$\Box$

\subsection{Nonexistence of radial solutions and maximum principle}

\subsubsection{Main results}
The following statements contribute to the nonexistence and triviality  for solutions to \eqref{zabka} and \eqref{zabkadiv}. 
 It can be treated as the problem overdermined by the condition $w(0)=0$.

\begin{tw}\label{twcopr}(Nonexistence and triviality for $C^1$ radial solutions) Let $B=B(0,R)\subseteq \mathbf{R}^n$, $R\in\mathbf{R}$, $n\ge 2$,  $w\! \in\ C^1(B)$ be the radial function such that
$w(0)=0$ and  $\Phi_p(\nabla w)\in W^{1,1}_{loc}(B, |x|^{-(n-1)}dx)$.
Moreover, suppose that one of the assumptions ($\mathcal{ND}_n$) or ($\mathcal{D}_n$) holds where: 
\begin{description}
\item[($\mathcal{ND}_n$)]
the set of conditions: $\mathcal{N}_{nd}:=
\{  \mathcal{A}1, \mathcal{A}2(a),\mathcal{A}3(a_l), \mathcal{A}4(a_l),  \mathcal{A}5(a_l)\}$, is satisfied  and $w$ is the solution to  \eqref{zabka}; 
\item[($\mathcal{D}_n$)]
the set of conditions $\mathcal{N}_{d}:= \{\mathcal{A}1,\mathcal{A}2(a),\mathcal{A}3(b_l), \mathcal{A}4(b_l),\mathcal{A}5(b_l)\}$ is satisfied and $w$ is the solution to  \eqref{zabkadiv}. 
\end{description}
Then  $w\equiv 0$. 
In particular, if $h(\cdot ,0,0)\neq 0$ on the set of positive measure in $(0,R)$, then such $w$ cannot  exist.
\end{tw}

Our next result is the variant of the maximum principle. As it asserts that under certain assumptions the solution has constant sign. It can be considered as generalization of a known theorem by Linqvist \cite{lindqvist}, when we deal with radial solutions.
Under some extra regularity assumptions, it also leads the nonexistence/triviality results.

\begin{tw}\label{twcole}(Maximum principle). Let $B=B(0,R)\subseteq \mathbf{R}^n$, $R\in\mathbf{R}$, $n\ge 2$, $w: B\rightarrow \mathbf{R}$ be the radial function such that  $w,\Phi_p(\nabla w)\in W^{1,1}(B\setminus B(0,r) )$ for any $0<r<R$   and
$$
w\equiv 0 \ {\rm on}\ \partial B\ \hbox{\rm in the sense of trace operator}.
$$ 
Moreover, suppose that one of the assumptions $(\mathcal{MND}_n )$ or 
$(\mathcal{MD}_n )$ holds where:
\begin{description}
\item[$(\mathcal{MND}_n )$]
the set of conditions: $\mathcal{M}_{nd}:=\left\{ \mathcal{A}1,\mathcal{A}2(b),\mathcal{A}3(a_r),\mathcal{A}4(a_r),\mathcal{A}5(a_r)\right\}$ is satisfied and $w$ is the solution to  \eqref{zabka}; 
\item[$(\mathcal{MD}_n )$]
the set of conditions $\mathcal{M}_{d}:=\left\{ 
\mathcal{A}1,\mathcal{A}2(b),\mathcal{A}3(b_r),\mathcal{A}4(b_r),\mathcal{A}5(b_r) \right\}$ is satisfied and $w$ is the solution to  \eqref{zabkadiv}. 
\end{description}
Then $w$ is of constant sign and monotone along the radii. Moreover,
\begin{displaymath}
\sup_{x\in B} |w(x)|=\limsup_{x\rightarrow 0}|w(x)|.
\end{displaymath}
If additionally $\left\{ w\in C(B)\right.$ and $\left.w(0)=0\right\}$ or
\begin{equation}\label{dopiskaa}
\limsup_{x\rightarrow 0} a(|x|)|\nabla w(x)|^p = 0
\end{equation}
then either $w\equiv 0$ or such $w$ cannot exist. 
\end{tw}

\begin{re}\rm 
The situation \eqref{dopiskaa} holds when for example $w\in C^1(B)$ and $a(\cdot)$ is bounded near $0$.
\end{re}

\noindent \textbf{Proof of Theorem \ref{twcopr}}: Let us denote  $w(x)=:u(|x|)$. 
According to Lemma \ref{lemreg}, we have 
$u\in W^{1,1}_{loc}((0,R))$. We will verify that $u$ fuilfills 
the assumptions in Theorem \ref{prawa}. For this, we observe at first that regularity assumption (a) is satisfied for $u$, by part 4 of Lemma \ref{lemreg}, and  boundary condition (b) also holds. Moreover, the regularity condition $w\in C^1(B)$ 
implies: 
$$
\nabla w (x)=u^{'}(|x|)\frac{x}{|x|}\ {\rm for}\ x\neq 0.$$
Therefore for any  $\theta\in S^{n-1}, r\in (0,R)$ we have $\nabla w(\theta r)= u^{'}(r)\theta$. As $\nabla w(\cdot)$ is continuous at zero, it implies that the limit $\lim_{r\to 0}u^{'}(r)\theta$ (which exists because $u^{'}\in C([0,R))$) is independent on $\theta$. This is possible only when $u^{'}(r)$  converges to zero at zero (and then $\nabla w(0)=0$). Thus $u(0)=u^{'}(0)=0$. The statement follows now from Theorem 2 applied to $u$ (part ii)), when we verify that in the case of ($\mathcal{ND}_n$) the assumption ($\mathcal{ND}$)
holds for $u$, while in the case of ($\mathcal{D}_n$) the assumption ($\mathcal{D}$)
is satisfied for $u$.
\hfill$\Box$

%\begin{re}\rm  
%For more general statement, we do not need to assume that  $w\in C^1(B)$. It is enough to assume instead that $w,\Phi_p(w) \in W^{1,1}_{loc} (B \setminus \{0\})$ and $\lim_{x\rightarrow 0}a(|x|)|\nabla w(x)|^p=0$. Therefore, nontrivial radial solution to the problem above may exist only when those conditions are violated. For example, if $a(\cdot)$ is continupus up to zero and  $a(0)=0$, then $\nabla w(\cdot)$ need not be bounded at zero. On the other hand, when $w\in C^1(B)$ is a radial nontrivial solution to the problem above, then $a(\cdot)$ cannot be bounded near zero.
%\end{re}

\bigskip
\noindent \textbf{Proof of Theorem \ref{twcole}}: 
Let $w(x)=: u(|x|)$, $u:(0,R)\rightarrow \mathbf{R}$. According to Lemma
\ref{lemreg}, we have $u,\Phi_p(u^{'})\in 
W^{1,1}_{loc}((0,R])$. Therefore $u$ fuilfills the assumptions (a) and (b)  in Theorem \ref{lewa} and moreover, $|u^{'}(|x|)|=|\nabla w(x)|$ a.e..
 Now it remins to note that when $(\mathcal{MND}_n )$ holds for $w$ then
$(\mathcal{MND})$ in Theorem \ref{lewa} holds for $u$, while if $(\mathcal{MD}_n )$ holds for $w$, then $(\mathcal{MD})$ in Theorem \ref{lewa} holds for $u$. 
% When $w\in C^1(B)$, then $w^{'}(0)=0$ and $\lim_{\epsilon\to 0}u^{'}(\epsilon)=0$. Consequently
%$\limsup_ {\epsilon\rightarrow 0} a(\epsilon)|u'(\epsilon)|^p = 0$ and we can apply last statement in Theorem \ref{lewa}.
\hfill$\Box$

\subsubsection{Examples within radial constraints}

Let us focus now on the  problem:
\begin{equation} \label{d1}
%\left\{ \begin{array}{l} 
 - |x|^\alpha\Delta_p w(x) + h(|x|,w(x), |\nabla w(x)|)=\phi(w(x))
 %\\
%\lim_{s\rightarrow 0}u(s)=0,\\
%\lim_{s\rightarrow 0} a(s)|u'(s)|^p=0 \end{array} \right.
\end{equation}
in $\mathcal{D}^{'}(B)$, where $B=B(0,1)\subseteq \mathbf{R}^n$ is the unit ball, $n\ge 2$, 
 $\Delta_pw(x)= \mathrm{div}(|\nabla w(x)|^{p-2} \nabla w(x))$ is the $p$-Laplacian, $1<p<\infty$ and $h$ satisfies:
 \begin{enumerate} 
 \item[(h)] 
 $(0,1)\times (\mathbf{R}\times [0,\infty)) \ni(\tau, (\lambda_0,\lambda_1))\mapsto h (\tau, \lambda_0,\lambda_1)$ is a Carath\'{e}odory function, i.e. it is measurable with respect to $\tau$ and continuous with respect to $(\lambda_0,\lambda_1)$, moreover
 $$|h(\t, \la_0 , \la_1)|\leq C\t^{\gamma} |\la_0|^l\la_1^{p-l-1}, \ {\rm where}\ 0<l<p,$$
  for every $\la_0\in\mathbb{R}, \la_1\in [0,\infty)$, almost every $\t\in(0,1)$, and where  $C>0$ is a given constant.
  \end{enumerate}

\subsubsection*{Nonexistence of radial solutions.}

The statement given below contributes to the nonexistence result 
from Theorem \ref{twcopr}.

\begin{tw} \label{twprzy}
Assume that
 \begin{enumerate}
\item[(a)] $1<p<\infty$, $0<l<p$, $0<\alpha<p$, $n\ge 2$, $n>\alpha(1-\frac{1}{p})+1$, $\gamma >-1$, $\gamma >\alpha -1-l$;
\item[(b)]  $\phi:\mathbf{R}\rightarrow \mathbf{R}$ is continuous, odd function such that $\tau \phi(\tau)>0$ a.e.;
\item[(c)] $h$ satisfies (h).
\end{enumerate}
Then there are no nontrivial radial solutions to  \eqref{d1} such that $w\in C^1(B)$, where $B=B(0,1)\subseteq\mathbf{R}^n$ is the unit ball,  $w(0)=0$ and\\ $\Phi_p(\nabla w)\in W^{1,1}_{loc}(B, 
|x|^{-(n-1)}dx)$.
\end{tw}

\noindent
\textbf{Proof}:\\
The assertion follows from Theorem \ref{twcopr}, under the assumptions 
$(\mathcal{ND}_n)$. We will show that the assumptions there are satisfied.
Clearly, $w$ is the solution to \eqref{zabka} where $a(\tau)=\tau^\alpha$. Moreover, the assumption $\mathcal{A}1$ is guranteed by (a), while 
$\mathcal{A}2(a)$ is the same as (b). When verifying  $\mathcal{A}3(a_l)$
we compute that $\delta_a(\tau)=C_{\alpha,p}\tau^{\alpha-1}$, where 
$C_{\alpha,p}:=(n-1)-\alpha(1-\frac{1}{p})>0$. Moreover,
\begin{eqnarray*}
a&\in& W^{1,1}_{loc}([0,1)) \Longleftrightarrow \alpha >0, \ 
\ \delta_a\in L^{1}_{loc}([0,1)) \Longleftrightarrow \alpha >0,\\
\ \delta_a^{-1/(p-1) }&\in & L^{1}_{loc}([0,1)) \Longleftrightarrow p>\alpha .
\end{eqnarray*}
Therefore \eqref{wyjazd}  in $\mathcal{A}3(a_l)$ is guaranteed by
 the assumption (a). 
To verify $\mathcal{A}4(a_l)$ we note that we have \eqref{a4} with 
$q(\tau)= C \tau^\gamma$, $\theta =1$, $C >0$. Moreover,  
$
q\in  L^{1}_{loc}([0,1)) \Longleftrightarrow \gamma >-1,
$ which is guaranteed by (a). 

When verifying $\mathcal{A}5(a_l)$ we recall that $q(\tau) = C_\gamma \tau^\gamma$, $v(\tau)=\delta_a(\tau)=C_{\alpha,p}\tau^{\alpha-1}$.
Hence
\begin{eqnarray*}
K &:=& K(0,r,q,v)= \left( \frac{p-l}{p} \right)^{\frac{p-l}{p}}\cdot \mathcal{L}^{l/p},\ {\rm where}\\
\mathcal{L} &:=&  \int_0^r \left( C_\gamma t^\gamma \right)^{p/l}
\left( C_{\alpha,p} t^{\alpha -1} \right)^{-(p-l)/l}
\left( \int_0^t   \left( C_{\alpha,p} \tau^{\alpha -1} \right)^{-1/(p-1)}
d\tau   \right)^{p-1} dt   \\
&\sim & \int_0^r t^\delta \left( \int_0^t \tau^\kappa d\tau    \right)^{p-1}dt,\  {for}\ 
\delta = \gamma \frac{p}{l} - (\alpha -1)\frac{p-l}{l},\ \ 
\kappa = -\frac{\alpha -1}{p-1}.
\end{eqnarray*}
Consequently
$$
K\sim \left[ \int_0^rt^{\delta +(\kappa +1)(p-1)}dt  \right]^{\frac{l}{p}}<\infty
\Longleftrightarrow \delta +(\kappa +1)(p-1) >-1 \Longleftrightarrow 
\gamma >\alpha -1-l,
$$
which holds because of (a). Therefore the statement follows.
\hfill $\Box$

\bigskip
\noindent
We will illustrate the above statement on the following example.

\begin{ex}[Sharpness of the assumption $\gamma >\alpha -l-1$]\rm~\\
Let us consider the function
$$ w(x):=|x|^s\ {\rm where}\ s>1.$$
An easy verification shows that 
$w\in  C^1(B), w(0)=0$.
 and 
\begin{eqnarray*}
\nabla w(x) &=& s|x|^{s-1}\frac{x}{|x|},\\
\Phi_p(\nabla w(x))_i &=& s^{p-1}\left( |x|^{(s-1)(p-1)}\frac{x_i}{|x|}\right)
=: s^{p-1}\left( |x|^{\kappa}\frac{x_i}{|x|}\right),\\
\\ \kappa &=& (s-1)(p-1),\\
\frac{\partial}{\partial x_i} \left(|x|^\kappa \frac{x_i}{|x|}  \right) &=& (\kappa-1)|x|^{\kappa -1}\frac{x_i^2}{|x|^2} + |x|^{\kappa-1},\\
\Delta_p w(x) &=& {\rm div}\left( \Phi_p(\nabla w(x))  \right) =s^{p-1}(\kappa -1+n)|x|^{\kappa -1},\\
-|x|^\alpha\Delta_p w(x) &=& - A |x|^{\kappa -1+\alpha},  \ 
A= s^{p-1}(\kappa -1 +n) >0.
 \end{eqnarray*}
In particular $\Phi_p(\nabla w(x)) \in 
 W^{1,1}_{loc}(B, |x|^{-(n-1)}dx)$.
We choose 
\begin{equation}\label{hallo}
h(\tau, \lambda_0,\lambda_1):= Ds^l \tau^\gamma |\lambda_0|^l|\lambda_1|^{p-l-1},%\ {\rm where}\ \gamma:= \alpha -l-2.
\end{equation}
where constant $D$ will be established later and verify that
\begin{eqnarray*}
h(|x|,w(x), |\nabla w(x)| ) =
Ds^{p-1}|x|^{\kappa +\gamma + l }.%=:B |x|^{\kappa-1+\alpha }.
%,\\
%\kappa = (s-1)(p-1)-1 .
%,\\\kappa_2 &=& \gamma +sl +(s-1)(p-l-1) = \kappa +l 
\end{eqnarray*}
%When we chose $|D|$ large enough 
Thus, for almost every $x\in B(0,1)$
$$
-|x|^\alpha\Delta_p w(x) + h(|x|,w(x),| \nabla w(x)| ) = s^{p-1}|x|^\kappa\left\{ -A|x|^{\alpha-1} +D|x|^{\gamma+l}\right\}.
$$
Right hand side above is of the form 
$\phi(w(x))$ where $\phi(\t )>0$ a.e. for $\t>0$ 
as in (b), 
if and only if
$$
\left\{ \gamma < \alpha -l-1 \ {\rm and}\ D\ge A\right\} \ {\rm or}\ 
\left\{ \gamma = \alpha -l-1 \ {\rm and}\ D> A\right\} 
$$
In particular, when the assumption $ \gamma > \alpha -l-1 $ does not hold, there are nontrivial solutions of \eqref{d1}, which satisfy the remaining assumptions in Theorem \ref{twprzy}. 
This shows sharpness of the assumption:  $\kappa -1+\alpha >0$.
 \end{ex}

\subsubsection*{Monotonicity property.}

Let us now consider Theorem \ref{twcole} within radial constraints.

\begin{tw}[Monotonicity]\label{monohomo}
Assume that  

\begin{enumerate}
\item[(a)]
$1<p<\infty$, $0<l<p$, $\alpha <p$, $n\ge 2, n> (1-\frac{1}{p})\alpha +1$, $\gamma = \alpha -1$;
\item[(b)]  $\phi:\mathbf{R}\rightarrow \mathbf{R}$ is continuous, odd function such that $\tau \phi(\tau)<0$ for a.e. $\tau$;  
\item[(c)] $h$ satisfies (h) with positive constant $C$ such that
$C\le \frac{X}{Y}$ where
\begin{eqnarray*}
X&:=& p^{1+\frac{l-1}{p}} (p-\alpha)^{(p-1)\frac{l}{p}}\left\{  
n-1-\alpha (1-\frac{1}{p})\right\} 
,\\
Y&:=& (p-1)^{(1-\frac{1}{p})(l+1)} .
\end{eqnarray*}
\end{enumerate}
Moreover, let
  $w: B\rightarrow \mathbf{R}$ be the radial function, where $B=B(0,1)\subseteq \mathbf{R}^n$ is the unit ball, 
     $w,\Phi_p(\nabla w)\in W^{1,1}(B\setminus B(0,r) )$ for any $0<r<R$,   
$$
w\equiv 0 \ {\rm on}\ \partial B\ \hbox{\rm in the sense of trace operator}
$$
and $w$ satisfies \eqref{d1}.

Then $w$ is of constant sign and monotone along the radii. Moreover,
\begin{displaymath}
\sup_{x\in B} |w(x)|=\limsup_{x\rightarrow 0}|w(x)|.
\end{displaymath}
If additionally
\begin{equation}\label{dopiskaaa}
\limsup_{x\rightarrow 0} |x|^\alpha|\nabla w(x)|^p = 0,
\end{equation}
then either $w\equiv 0$ or such $w$ cannot exist. 
\end{tw}

\noindent \textbf{Proof}:\\ 
The statement follows from Theorem \ref{twcole} and we have to verify the assumptions therein.
Clearly, $w$ satisfies \eqref{zabka} with $a(\tau)=\tau^\alpha$, after we note that $|\langle \nabla w(x),\frac{x}{|x|}\rangle =|\nabla w(x)|$. Moreover, the regularity assumptions also hold for $w$. We have to confirm that
 the set of conditions 
$\mathcal{MND}_n$ is satisfied when 
$q(\tau)=C\tau^\gamma$, $h$ is as in $(h)$ and $R=1$.
We only verify part $\mathcal{A}5(a_r)$, laving the remaining verifications to the reader. 

We have $$\delta_a(\tau) = C_{\alpha,p}\tau^{\alpha-1}\ {\rm where 
}\ C_{\alpha,p}:=(n-1)-\alpha(1-\frac{1}{p})>0,$$ so that
$K:= K(0,1, C\tau^\gamma, C_{\alpha,p}\tau^{\alpha-1})=$
\begin{eqnarray*}
A
\left\{   
\int_0^1 (Ct^{\alpha -1})^{p/l} (C_{\alpha,p}\tau^{\alpha-1})^{-(p-l)/l}
\left( \int_0^t (C_{\alpha,p}\tau^{\alpha-1})^{-1/(p-1)}d\tau \right)^{p-1} dt
\right\}^{l/p},
\end{eqnarray*}
where $A=\left( \frac{p-1}{p} \right)^{ \frac{p-1}{p} }$.
Thus, to have $K<\infty$, we have to require that $(1-\alpha)/(p-1)>-1$, equivalently $\alpha <p$. Then
\begin{eqnarray*}
K= B \left( \int_0^1 t^\kappa dt\right)^{l/p} = \frac{B}{p^{l/p}} \stackrel{(c)}{\le} 1,
\end{eqnarray*}
where $B=\left( \frac{p-l}{p} \right)^{ \frac{p-l}{p} }C(C_{\alpha,p})^{-1}(\frac{p-1}{p-\alpha})^{(p-1)l/p}$,
\begin{eqnarray*}
\kappa &=& (\alpha -1) \frac{p}{l}-\frac{(\alpha-1)(p-l)}{l}+ (-\frac{\alpha -1}{p-1} +1)(p-1) 
%&=& \frac{p}{l}(\gamma -\alpha +1) +\alpha -1 -\alpha +p \\
=  p-1\  (>-1). 
\end{eqnarray*}
\hfill$\Box$

\bigskip
\noindent
We will illustrate the above result within the restricted class of functions.

\begin{ex}[Confirmation within the restricted class]
Let $n=2$ and 
\begin{equation}\label{forma}
w(x):= 1-|x|^s, s>0,
\end{equation}
where parameter  $s$ will be specified later. Obviously, $w$ decreases 
up to $0$ achieved on $\partial B$ and it is nontrivial.
%Then $w\in W^{1,1}(B)\cap C^1(B)$, $\Phi_p(\nabla w)\in W^{1,1}(B)$ and $|x|^\alpha$ is bounded near zero, but $\tilde{w}$ is nontrivial.
%We will expalin that $w$ cannot fuilfill the assusmptions 
%(a)-(c) in Theorem \ref{monohomo}. 
We compute that:
\begin{eqnarray*}
\nabla w (x)=-s|x|^{s-1}\frac{x}{|x|},\ \Phi_p(\nabla w(x))=-
s^{p-1}|x|^{\kappa}\frac{x}{|x|},\ 
\kappa = (s-1)(p-1);\\
\frac{\partial}{\partial x_j} \left(\Phi_p(\nabla w(x))_i\right) = - s^{p-1}|x|^{\kappa -1}\left\{ (\kappa -1) \frac{x_ix_j}{|x|^2} +\delta_{ij}\right\},\\
\Delta_pw(x)=- s^{p-1}\left( (\kappa + 1)|x|^{\kappa -1}\right),\\
-|x|^\alpha \Delta_pw(x) =s^{p-1}(\kappa + 1) |x|^{\kappa +\alpha -1}.
\end{eqnarray*}
In particular such $w$ satisfies the  regu,arity assumptions in Theorem 
\ref{monohomo}. 

With the same $h$ as in \eqref{hallo}, where $\gamma =\alpha -1$, we get 
\begin{eqnarray*}
h(|x|,w(x), |\nabla w(x)| ) =
Ds^{p-1} (1-|x|^s)^l|x|^{ \kappa+ \alpha -1          -l(s-1) },\\
-|x|^\alpha \Delta_pw(x) + h(|x|,w(x), |\nabla w(x)| )=~~~~~~~~~~~~~~~~~~~~~~~~~~~~\\ =
s^{p-1} |x|^{\kappa +\alpha -1-l(s-1)}\left\{ (\kappa + 1)|x|^{l(s-1)} +D(1-|x|^s)^l \right\}.
\end{eqnarray*}
As $|x|=(1-w(x))^{1/s}$, this is in the form 
$\phi(w(x))$ when 
$$
\phi(w) = s^{p-1}(1-w)^{\frac{\kappa +\alpha -1-l(s-1)}{s}  }\left\{
(\kappa +1) (1-w)^{\frac{l(s-1)}{s} } +Dw^l
\right\} 
$$
on $[0,1]$.

%Let us start with the verification of condition (h).
Note that by our assumptions we need $\phi <0$ a.e..
When verifying the sign of $\phi$ near $0$ and $1$, we note that it is possible only when 
\begin{equation}\label{z1}
\kappa \le -1   \Leftrightarrow s\le  \frac{p-2}{p-1} (<1), \ \  D<0.
\end{equation}
In that case the function
$$
v(w) := (\kappa+1) + Dw^l (1-w)^{l(1-s)/s}
$$
is negative a.e. $[0,1]$ and so if $\phi$.

Obviously, $\phi$ is continuous at $0$, while for continuoty of
$\phi$ at $1$ we require
\begin{equation}\label{z2}
s+\alpha -1\ge 0\Leftrightarrow s\ge 1-\alpha.
\end{equation}

Linking \eqref{z2} with \eqref{z1} and the condition $\alpha <p$ from Theorem \ref{monohomo}, (a), we obtain
\begin{equation}\label{z3}
p>\alpha \ge \frac{1}{p-1}. 
\end{equation}  
  In our case of $n=2$, so the condition $n> (1-\frac{1}{p})\alpha +1$ in (a) reads 
\begin{equation}\label{z4}
\alpha <\frac{p}{p-1}  .
\end{equation}  
Obviously (h) will be satisfied with $C:=(-D)$, when choosing  $D$ not to small. 

Taking into account the involved conditions we then get
\begin{eqnarray*}
\mathcal{G}:= 
\alpha + p(s-1)\stackrel{\eqref{z1}}{\le} \alpha +p\left( \frac{p-2}{p-1}-1  \right) =\alpha -\frac{p}{p-1}\stackrel{\eqref{z4}}{<}   0, 
\end{eqnarray*}
while the condition \eqref{dopiskaa} is equivalent to $
\mathcal{G} >0 
%\Longleftrightarrow s-1>\frac{-\alpha}{p}
%\Longleftrightarrow \kappa > (-\alpha )\frac{p-1}{p}
%\Longleftrightarrow \kappa +1 > (-\alpha ) \frac{p-1}{p} +1 =\mathcal{T}
.$

Thus,  functions  like \eqref{forma} are example illustrations of  Theorem \ref{monohomo} within nontrivial function. On the other hand, 
we confirm  that there are no such  functions
 satisfying additionally the condition \eqref{dopiskaaa}.

%so when the above holds, then functions like \eqref{forma}
%cannot verify our assumptions (a)-(c)

\end{ex}

\section{Final remarks}

We end our discussion with the following remarks.

\begin{re}\rm
In all our presented problems dealing with nonexistence/triviality, the growth of $\phi$ does not play any role.
\end{re} 

\begin{re}\rm
Further development of Opial-type inequalities leds to generalizations of our results, where one can consider more general class of nonlinearities $h(\cdot,\cdot,\cdot)$. Moreover, it is possible also to obtain generalizations of our results within the class of $A$-harmonic problems
like
\begin{displaymath}
 - \mathrm{div}(a(|x|) A(\nabla w)) + h(|x|,w(x),\langle\nabla w(x),\frac{x}{|x|}\rangle)=\phi(w(x)), 
\end{displaymath} 
where $A:\mathbf{R}^n\rightarrow \mathbf{R}^n$ is given function.
\end{re}

%\newpage


\begin{thebibliography}{99}
\bibitem{AdKa2009} T. Adamowicz, A. Ka\l{}amajska, {\it On a variant of the maximum principle involving radial $p$-Laplacian with applications to nonlinear eigenvalue problems and nonexistence results},  Topol. Methods Nonlinear Anal. {\bf 34},  (2009) 1–20.
\bibitem{AdKa2010} T. Adamowicz, A. Ka\l{}amajska {\it Maximum principles and nonexistence results for radial solutions to equations involving \mbox{$p$-Laplacian}},  Math. Methods Appl. Sci. . {\bf 33}(13) (2010), 1618-1627.
%\bibitem{AgarD}  R. P. Agarwal, {\it Difference Equations and Inequalities}, Marcel Dekker Inc., New York, 1992.
%\bibitem{AgarPa} R. P. Agarwal, P. Y. Pang {\it Opial Inequalities with Applications in Differential and Difference Equations}, 1995.
%\bibitem{AgarPa95} R. P. Agarwal, P. Y. Pang, {\it Sharp Opial-type inequalities in two variables} Appl. Anal. 56 (1995), nr 3-4, 227–242.
%\bibitem{AgarPa97} R. P. Agarwal, P. Y. Pang, {\it Opial-type inequalities involving higher order partial derivatives of two functions}, General inequalities 7 (Oberwolfach, 1995), 157–178, Internat. Ser. Numer. Math., 123, Birkh\"{a}user, Basel, 1997.
%\bibitem{AgaShe}  R. P. Agarwal, Q. Sheng, {\it Sharp integral inequalities in $n$-independent variables}, Nonlinear Anal. 26 (1996), nr 2, 179–210.
%\bibitem{Almenar}  P. Almenar, L. Jódar {\it Asymptotic behaviour of the solutions of second order functional differential equations} Comput. Math. Appl. 62 (2011), nr 1, 297–309.
\bibitem{ArDiTe} D. Arcoya, J. Diaz and L. Tello, {\it S-shaped bifurcation branch in a quasilinear multivalued model arising in climatology},  J. Differential Equations {\bf 150} (1998), 215--225.
%\bibitem{Ba} J. M. Ball, {\it Some Open Problems in Elasticity, Geometry, Mechanics and Dynamics}, Springer, New York, 2002.
%\bibitem{BaFaHo} J. Batt, W. Faltenbacher, E. Horst, {\it Stationary spherically symmetric models in stellar dynamics}, Archive for Rational Mechanics and Analysis 93 (1986), 159–183.
%\bibitem{Bee}  P. R. Beesack, {\it Elementary proofs of some Opial-type integral inequalities}, J.~d'Analyse Mathematique 36 (1979), 1-14.
\bibitem{BeeDas} P. R. Beesack, K. M. Das {\it Extensions of Opial inequality}, Pacific J. Math. 26 (1968), 215-232.
\bibitem{BDFP} V.  Benci,  P.  D'Avenia,  D.  Fortunato  and  L.  Pisani, {\it Solitons  in  several  spacedimensions:  Derrick's problem and infinitely many solutions}, Arch. Rat. Mech. Anal. 154 (2000), 297--324.
%\bibitem{BeFoPi} V. Benci, D. Fortunato and L. Pisani, {\it Soliton like solutions of a Lorentz invariant equation in dimension 3}, Rev. Math. Phys. 10 (1998), 315--344.
%\bibitem{Bird} R. B. Bird, W. E. Stewart, E. N. Lightfoot {\it Transport Phenomena}, Second Edition, John Wiley \& Sons, Inc.
%\bibitem{BoydWong} D. W. Boyd, J. S. W. Wong, {\it An extension of Opial's inequality}, J. Math. Anal. Appl. 19 (1967), 100-102.
%\bibitem{ra1} F. Brock, {\it Radial symmetry for nonnegative solutions of semilinear elliptic equations involving the $p-$Laplacian}, Proceedings of the Conference Calculus of Variations, Applications and Computations, Pont-\`a-Mousson, 1997.
%\bibitem{te1} J. F. Brothers, W. P. Ziemer, {\it Minimal rearrangements of Sobolev functions}, J.~Reine Angew. Math. 384 (1988), 153–179.
%\bibitem{Calvert} J. Calvert, {\it Some generalizations of Opial's inequality}, Proc. Amer. Math. Soc. 18 (1967), 72-75.
\bibitem{bfg} M.F. Bidaut-V\'eron, M. García-Huidobro, L. Véron, {\em 
Estimates of solutions of elliptic equations with a source reaction term involving the product of the function and its gradient,}
Duke Math. J. {\bf 168}(8) (2019),  1487–1537. 
\bibitem{CaNa} A. Callegari, A. Nachman {\it A nonlinear singular boundary value problem in the theory of pseudoplastic fluids}, SIAM J. Appl. Math. 38(1980), 275--281.
%\bibitem{te2} A. Cianchi, A. Ferone, {\it On symmetric functionals of the gradient having symmetric equidistributed minimizers}, J. Math. Anal. 38 (2006), 279–308.% SIAM J.

\bibitem{cm} G. Caristi, E. Mitidieri, {\em 
Nonexistence of positive solutions of quasilinear equations,}
Adv. Differential Equations 2 (1997), no. 3, 319–359. 



%\bibitem{ra2} L. Damascelli, F. Pacella, {\it Monotonicity and symmetry of solutions of $p-$Laplace equations, 1 $< p <$ 2, via the moving plane method}, Ann. Scuola Norm. Sup. Pisa Cl. Sci. 26 (1998), 689–707.
\bibitem{Diaz} J. I. Diaz {\it Nonlinear Partial Differential Equations and Free Boundaries}, Vol. I, Elliptic Equations, Res. Notes Math., vol. 106, Pitman Publ. Ltd., Boston, London, Melbourne, 1985.
\bibitem{Die} L. Diening, A. Prohl and M.  Růžička, {\it Semi-implicit Euler scheme for generalized Newtonian fluids}, SIAM J. Numer. Anal. 44 (2006), 1172--1190.
%\bibitem{ra3} J. Dolbeault, P. Felmer, R. Monneau, {\it Symmetry and nonuniformly elliptic operators}, Differential Integral Equations 18 (2005), 141–154.
\bibitem{DOIw} L. D'Onofrio and T. Iwaniec, {\it The $p$-harmonic transform beyond its natural domain of definition},  Indiana Univ. Math. J.  53 (2004), 683--718.
\bibitem{Dra1} P. Drábek, {\it The $p$-Laplacian -- Mascote of Nonlinear Analysis}, Acta Math. Univ. Comenianae, Vol. LXXVI (Proceedings of Equadiff 11)
1 (2007), 85--98. 
% \bibitem{Dra} P. Drábek, M. García-Huidobro, R. Manásevich, {\it Positive solutions for a class of equations with a $p$-Laplace like operator and weights}, Nonlinear Anal. 71 (2009), nr 3-4, 1281-1300.
%\bibitem{wir} H. Dym, H. McKean, {\it Fourier series and integrals}, Academic press, 1985. 
%\bibitem{Evans} L. C. Evans {\it Równania różniczkowe cząstkowe}, Wydawnictwo Naukowe PWN. Warszawa 2002.
%\bibitem{Dra01} M. García-Huidobro, A. Kufner, R. Manásevich, C. Yarur, {\it Radial solutions for a~quasilinear equation via Hardy inequalities}, Differential and Integral Equations 6 (2001) 1517–1540.
\bibitem{DraKuNi} P. Drábek, A. Kufner, F. Nicolosi, {\it Quasilinear elliptic equations with degenerations and singularities} De Gruyter Series in Nonlinear Analysis and Applications, 5. Walter de Gruyter \& Co., Berlin, 1997. 
%\bibitem{Ev} L. C. Evans {\it Partial Differential Equation} Graduate studies in mathematics, vol 19, American Mathematical Society, 1998.
{\bibitem{filipucci1} R. Filippucci, {\em Nonexistence of positive weak solutions of elliptic inequalities,} Nonlinear Anal. 70 (2009), 2903--2916.}
\bibitem{FoOrPi} D. Fortunato, L. Orsina and L. Pisani, {\it Born--Infeld type equations for electrostatic fields}, J. Math. Phys. 43 (2002), 5698--5706.
%\bibitem{GiNi} B. Gidas, W. M. Ni, L. Nirenberg, {\it Symmetry and related properties via the maximum principle}, Commun. Math. Phys. 68 (1979), 209–243.
%\bibitem{te3} A. Gladkov, N. Slepchenkov, {\it Entire solutions of quasilinear elliptic equations}, Nonlinear Anal. 66 (2007), 750–775.
%\bibitem{YangHwang90} T. S. Hwang, G. S. Yang, {\it On integral inequalities related to Opial's inequality}, Tamkang J. Math. 21 (1990), 177-183.
%\bibitem{LiKa} A. Ka\l{}amajska, K. Lira {\it Maximum modulus principles for radial solutions of quasilinear and fully nonlinear singular P.D.E’s},  Bull. Belg. Math. Soc. 14, (2007) 157-176.
\bibitem{KaSt} A. Ka\l{}amajska, A. Stryjek, {\it On maximum principles in the class of oscillating functions}, Aequationes Math. 69 (2005), 201-211.
\bibitem{KaYaYo} N. Kawano, E. Yanagida and S. Yotsutani, {\it Structure theorems for positive radial solutions to} div$(|Du|^m Du) + K(|x|)u^q= 0$ in $\mathbb{R}^n$, J. Math. Soc. Japan 45 (1993), 719--742.
%\bibitem{Ka} B. Kawohl, {\it Rearrangements and convexity of level sets in PDE}, Lecture Notes in Mathematics, vol. 1150, Springer--Verlag, Berlin--Heidelberg, 1985.
%\bibitem{Bear} R. Laister, R. Beardmore, {\it Transversality and separation of zeros in second order differential equations},  Proc. Amer. Math. Soc. 131 (2003), no. 1, 209–218.
%\bibitem{Lasota}  A. Lasota, {\it A discrete boundary value problem}, Ann. Polon. Math. 20 (1968), 183-190.

%\bibitem{ra4} Y. Li, W. M. Ni, {\it Radial symmetry of positive solutions of nonlinear elliptic equations in\ $\mathbb{R}^n$}, Comm. Partial Differential Equations 18 (1993), 1043–1054.
%\bibitem{Mat} T. Matukuma, {\it The Cosmos}, Iwanami Shoten, 1938.
%\bibitem{Maz} V. G. Maz’ya, {\it Sobolev Spaces}, Springer-Verlag, 1985.
% \bibitem{MitriPeca} D. S. Mitrinović, J. E. Pe\u{c}arić, {\it Generalizations of two inequalities of Godunova and Levin}, Bull. Polish Acad. Sci. Math. 36 (1988), 645-648.
% \bibitem{Necaev}  I. D. Ne\u{c}aev, {\it Integral inequalities with gradients and derivatives}, Soviet Math. Dokl. 22 (1973), 1184-1187.
% \bibitem{Olech} C. Olech, {\it A simple proof of a certain result of Z. Opial}, Annales Polonici Mathematici nr 1, tom 8 (1960), 61-63.
% \bibitem{Opia} Z. Opial, {\it Sur une inégalité de C. de la Vallée Poussin dans la théorie de l'équation différentielle linéaire du second ordre}, Ann. Polon. Math. 6 (1959/1960), 87-91.
\bibitem{kuf-opic} A. Kufner, B. Opic, {\em How to define reasonably weighted sobolev spaces,} Comment. Math. Univ. Carolin., {\bf 25}(3) (1984), 537--554.
\bibitem{Lee}  C. S. Lee, {\it On some generalization of inequalities of Opial, Yang and Shum}, Canad. Math. Bull. 23 (1980), 71-80.
\bibitem{lindqvist} P. Lindqvist, {\em On the equation ${\rm div}\,(|\nabla u|^{p-2}\nabla
              u)+\lambda|u|^{p-2}u=0$,} Proc. Amer. Math. Soc. {\bf {109}}(1) (1990), 157--164.
\bibitem{pohmi_99_b} E. Mitidieri and S. Pohozaev,  {\em A priori estimates and blow-up of solutions to nonlinear partial differential equations and inequalities. Transl. from the Russian},
Proc. Steklov Inst. Math. {\bf 234} (2001), 1--362. (translated
from Tr. Mat. Inst. Steklova 234 (2001), 1--383.)
              
 \bibitem{Opial} Z. Opial, {\it Sur une inégalité}, Ann. Polon. Math. 8 (1960), 29-32.
% \bibitem{Pach} B. G. Pachpatte, {\it On Opial-type integral inequalities}, J. Math. Anal. Appl. 120 (1986), 547-556.
 \bibitem{PachD87} B. G. Pachpatte, {\it A note on Opial and Wirtinger type discrete inequalities}, J.~Math. Anal. Appl. 127 (1987), 470-474.
% \bibitem{Pach89}  B. G. Pachpatte, {\it On certain two dimensional integral inequalities}, Chinese J.~Math. 17 (1989), 273-279.
% \bibitem{PachD90} B. G. Pachpatte, {\it On Opial type discrete inequalities}, Anal. sti. Univ. "Al. I. Cuza" din Iasi 36 (1990), 237-240.
% \bibitem{QiZhe} Z. Qi, {\it Further Generalizations of Opial's Inequality}, Acta Mathematica Sinica, New Series, tom 1, nr 3 (1985) 196-200.
%\bibitem{Ru} M.  Ru\u02c7zi\u02c7cka, {\it Electrorheological  Fluids:  Modeling  and  Mathematical  Theory, Lecture Notes in Mathematics}, vol. 1748, Springer, Berlin, Germany, 2000.
% \bibitem{Rozanova} G. I. Rozanova, {\it Ob odnom integral'nom neravenstve, svjazannom s~neravenstvom Polia}, Izvcstija Vyss. Ucebn. Zaved. Mat. 125 (1975), 75-80.
% \bibitem{Saker} S. H. Saker, {\it Lyapunov's type inequalities for fourth-order differential equations} Abstr. Appl. Anal., 2012.
%\bibitem{ra5} J. Serrin, H. Zou, {\it Symmetry of ground states of quasilinear elliptic equations}, Arch. Rational Mech. Anal. 148 (1999), 265–290.
\bibitem{sin91} G. J. Sinnamon, {\em 
Weighted Hardy and Opial-type inequalities,}
J. Math. Anal. Appl. {\bf 160}(2) (1991), 434--445. 
 \bibitem{Shum}  D. T. Shum, {\it On a class of new inequalities}, Trans. Amer. Math. Soc. 204 (1975), 299-341.
% Math. Z. 4 (1919), 139-151.
% % \bibitem{Willet} D. Willett, {\it The existence-uniqueness theorem for an n-th order linear ordinary differential equation}, Amer. Math. Monthly 75 (1968), 174-178.
% \bibitem{Wong} J. S. W. Wong, {\it A discrete analogue of Opial's inequality}, Canad. Math. Bull. 10 (1967), 115-118.
% \bibitem{Yang66} G. S. Yang, {\it On a certain result of Z. Opial}, Proc. Japan Acad. 42 (1966), 78-83.
% \bibitem{Yang87}  G. S. Yang, {\it A note on an inequality similar to Opial inequality}, Tamkang J. Math. 18 (1987), 101-104.
\bibitem{WuZhLi} Z. Wu, J. Zhao and H. Li, {\it Nonlinear Diffusion Equations}, World Scientific, Singapore, 2001.
\bibitem{Sze1} G. Szeg\"o, {\it Orthogonal polynomials, American Mathematical Society Colloquium Publications}, 1939.
%\bibitem{Sze2} G. Szeg\''o, {\it \"Uber Orthogonalsysteme von Polynomen}, Ibid., vol. 4 (1919), 139--151.	
%\bibitem{Maz} V. Maz’ya {\it Sobolev Spaces with Applications to Elliptic Partial Differential Equations}, 1985
\end{thebibliography}
\end{document}